\documentclass[10pt]{article}
\usepackage{amsmath}
\usepackage{amsfonts}
\usepackage{latexsym}
\usepackage{amssymb}
\usepackage[dvips]{graphics}
\newcommand{\qed}{\begin{flushright} $\Box$\end{flushright}}

\newcommand{\Z}{{\mathbf Z}}

\newcommand{\R}{{\mathbf R}}

\newcommand{\cat}{{\rm {cat }}}

\newcommand{\id}{{\rm {id }}}

\newcommand{\tc}{{\rm\bf {TC}}}
\newcommand{\proofend}{$\Box$\bigskip}

\newcommand{\RP}{\mathbf {RP}}
\newcommand{\CP}{\mathbf {CP}}
\newtheorem{theorem}{Theorem}
\newtheorem{prop}[theorem]{Proposition}
\newtheorem{lemma}[theorem]{Lemma}
\newtheorem{corollary}[theorem]{Corollary}
\newtheorem{definition}{Definition}
\newtheorem{remark}{Remark}
\newenvironment{proof}{\par\noindent{\bf Proof.} }{\par\noindent}

\title{Topological Robotics: Motion  Planning in Projective Spaces}
\author{Michael Farber\footnote{Partially supported by
the US - Israel Binational Science Foundation.}\, \footnote{Partially supported by the Herman Minkowski Center for Geometry}\, , \, Serge
Tabachnikov$^\ast$ and Sergey Yuzvinsky$^\ast$}
\date{October 2, 2002}

\begin{document}
\maketitle

\section{Introduction}

In this paper we study one of the most elementary problems of the topological robotics: rotation of a line,
which is fixed by a revolving joint at a base point: 
one wants to bring the line from its initial position $A$ to a final position $B$ by a
continuous motion in the space. The final goal is to construct a motion planning
algorithm which will perform this task once the initial
position $A$ and the final position $B$ are presented. This problem becomes hard when the dimension of the space 
is large.

Any such motion planning algorithm must have instabilities, i.e. 
the motion of the system will be discontinuous as a function of $A$ and $B$. 
These instabilities are caused by topological reasons.
A general approach to study instabilities of robot motion was suggested recently in papers \cite{F1, F2}.
With any path-connected topological 
space $X$ one associates in \cite{F1, F2} a number $\tc(X)$, called {\it the topological complexity of $X$}. 
This number is of fundamental importance for the 
motion planning problem: $\tc(X)$ determines character of instabilities which have all  
motion planning algorithms in $X$.

The motion planning problem of moving a line in $\R^{n+1}$ reduces to a topological problem of 
calculating the topological complexity of the real projective space $\tc(\RP^n)$, which we tackle in this paper. 
We compute the number $\tc(\RP^n)$ for all $n\leq 23$ (see the table in section \ref{mres}). 
This will probably be useful for applications.

Surprisingly, the problem of finding a general formula for $\tc(\RP^n)$ turns out to be quit difficult. 
One of the main results
of this paper claims that for $n\not=1, 3, 7$, the number 
$\tc(\RP^n)$ coincides with the smallest integer $k$ such that the projective space $\RP^n$ admits an immersion into $\R^{k-1}$.
This means that our problem of finding $\tc(\RP^n)$ is equivalent to the immersion problem for real projective spaces.
The latter is a famous, classical, well researched, topological problem, see \cite{Da, Ja1} for a survey. 
A general answer to the immersion problem $\RP^n\to \R^{k-1}$ is not known, but there are many 
important results in the literature.

For the remarkable special values $n=1, 3, 7$ the answer is very simple: $\tc(\RP^n) = n+1$. 

We show that for $n\leq 8$ explicit motion planning algorithms for lines in $\R^n$ can be build by 
using multiplication of the complex numbers, the quaternions, and the Cayley numbers.

\section{The motion planning problem}\label{section1}
\label{definition}

In this section we recall some definitions and results from \cite{F1, F2}. 

Let $X$ be a topological space $X$, thought of as the configuration
space of a mechanical system. 
We will always assume that $X$ is
a connected CW complex.
Given two points $A, B \in X$,
one wants to connect them by a path in $X$; this path represents a continuous motion of the system from
one configuration to the other. 
A solution to this motion
planning problem is a rule (algorithm) that takes $(A,B)
\in X \times X$ as an input and produces a path from $A$ to $B$ as an
output. Let $PX$ denote the space of all continuous paths $\gamma: [0,1] \to X$,
equipped with the compact-open topology, and
let $\pi: PX \to X \times X$ be the map assigning the end points to a
path: $\pi (\gamma) = (\gamma(0),
\gamma(1))$. The map $\pi$ is a fibration whose fiber is the
based loop space $\Omega X$.
The motion planning problem consists of finding a section $s$ of this
fibration.

The section $s$ cannot be continuous, unless $X$ is
contractible, see \cite{F1}. One defines $\tc(X)$, the {\it
topological complexity of} $X$, as the smallest
number $k$ such that $X \times X$ can be
covered by open sets $U_1, \dots, U_k$, so that for every $i=1,
\dots, k$ there exists a continuous section $s_i:
U_i \to PX, \pi \circ s_i = \id$.

According to \cite{F2}, {\it a motion planner} in 
$X$ is defined by finitely many subsets $F_1, \dots, F_k$ $\subset X\times X$ and continuous maps
$s_i: F_i\to PX$, where $i=1, \dots, k$, such that:
\begin{enumerate}
\item[(a)] the sets $F_1, \dots, F_k$ are pairwise disjoint $F_i\cap F_j=\emptyset$, $i\not= j$, and cover $X\times X$;
\item[(b)] $\pi\circ s_i =1_{F_i}$ for any $i=1, \dots, k$;
\item[(c)] each $F_i$ is an ENR.
\end{enumerate}

The subsets $F_i$ are {\it local domains} of the motion planner; the maps $s_i$ are {\it local rules}.
Any motion planner determines a {\it motion planning algorithm}:
given a pair $(A, B)$ of initial - final configurations, 
we determine the index $i\in \{1, 2, \dots, k\}$,
such that the local domain $F_i$ contains $(A, B)$; then we apply the local rule $s_i$ and 
produce the path $s_i(A,B)\in PX$ in $X$ as an output.

{\it The order of instability} of a motion planner is defined as the largest $r$ such that the closures of some $r$
among the local domains $F_1, \dots, F_k$ have a non-empty intersection: 
\begin{eqnarray*}
\bar F_{i_1}\cap \bar F_{i_2}\cap \dots \cap \bar F_{i_r}\, \not=\, \emptyset
\quad\mbox{where}\quad 1\leq i_1<i_2<\dots< i_r\leq k.
\end{eqnarray*} 
The order of instability describes the character of discontinuity of the motion planning algorithm
determined by the motion planner. 

In \cite{F2} it is shown that:
{\it the minimal integer $k$, such that a smooth manifold $X$ admits a motion planner
with $k$ local rules, 
equals $\tc(X)$.  Moreover,
the minimal integer $r>0$, such that $X$ admits a motion planner with order of instability $r$, equals $\tc(X)$.}

This explains importance of knowing the number $\tc(X)$ while solving practical motion planning problems.

Let us mention some other results from \cite {F1}, which will be used later in this paper.
$\tc(X)$ depends only on the homotopy type of $X$.
One has
\begin{eqnarray}\label{tab}
{\rm cat} (X) \leq \tc(X) \leq 2\ {\rm cat} (X) -1,
\end{eqnarray}
where ${\rm cat} (X)$ is the Lusternik-Schnirelmann 
category of $X$.
If $X$ is
\hfill\break
$r$-connected then
\begin{eqnarray}\label{stam}
\tc(X) < \frac{2\ {\rm dim}(X)+1}{r+1} + 1.
\end{eqnarray}

Next result provides a lower bound for $\tc(X)$ in terms of the
cohomology ring  $H^\ast(X)$ with coefficients in a
field. The tensor product $H^\ast(X)\otimes H^\ast(X)$ is also a graded ring
with the multiplication
$$
(u_1\otimes v_1)\cdot (u_2\otimes v_2) = (-1)^{|v_1|\cdot |u_2|}\,
u_1u_2\otimes v_1v_2
$$
where $|v_1|$ and $|u_2|$ are the degrees of the cohomology classes
$v_1$ and $u_2$. The cohomology multiplication
$H^\ast(X)\otimes H^\ast(X)\to H^\ast(X)\label{prod}$ is a ring
homomorphism. Let  $I \subset H^\ast(X)\otimes
H^\ast(X)$ be the kernel of this homomorphism. The ideal $I$ is
called {\it the ideal of zero-divisors} of
$H^\ast(X)$. The {\it zero-divisors-cup-length} is the length of the
longest nontrivial product in the ideal of
zero-divisors. It is shown in \cite{F1} that {\it the topological complexity $\tc(X)$ is greater
than the zero-divisors-cup-length of
$H^\ast(X)$.}

In this estimate, one can equally well use extraordinary cohomology
theories instead of the usual cohomology.

The topological complexity $\tc(X)$, as well as 
the Lusternik-Schnirelmann category $\cat(X)$, are particular cases of the notion of {\it Schwarz
genus} (also known as {\it sectional category}) of a fibration; it was introduced and 
thoroughly studied by A.Schwarz in \cite{Sch}, see also \cite{Ja2, Ja3}. 

\section{Motion planning in simply connected symplectic manifolds}

In this section we find the topological complexity of motion planning in any closed simply connected
symplectic manifold:

\begin{theorem}\label{complex}
Let $X$ be a closed $2n$-dimensional simply connected symplectic manifold.
Then $\tc(X)=2n+1$.
\end{theorem}

\begin{proof} Since $X$ is simply connected, the inequality  $\tc(X)
\leq 2n+1$ follows from (\ref{stam}). On the other hand, 
let $\omega$ be the symplectic 2-form on $X$.
If $[\omega] \in H^2(X;\R)$ denotes the respective cohomology
class then $[\omega]^n \neq 0$. Hence 
$[\omega]\otimes 1 - 1 \otimes [\omega]\in H^\ast(X;\R)\otimes H^\ast(X;\R)$
is a zero-divisor (see the definition in \S \ref{section1}) whose $2n$-th power
$$([\omega]\otimes 1 - 1 \otimes [\omega])^{2n}$$ 
does not vanish  
since it contains the term
$$(-1)^n\binom {2n} n [\omega]^n \otimes [\omega]^n.$$ 
The opposite
inequality $\tc(X) \geq 2n+1$ now follows from Theorem 7 of \cite{F1}. \proofend
\end{proof}

This result applies, in particular, to the complex projective space $\CP^n$:
\begin{corollary}\label{complex1} One has, $\tc({\bf CP}^n)=2n+1.$
\end{corollary}

Theorem \ref{complex} also applies to other smooth simply connected algebraic varieties
such as complex Grassmannians and flag varieties.

\begin{remark} {\rm A similar result, ${\rm cat}(X) = n+1$, for any 
closed simply connected symplectic manifold $X$,
is well
known \cite{Be1,Ja2}.}
\end{remark}

\section{Motion planning in the real projective space: the initial discussion}\label{initial}

In this section we start studying the problem of computing the topological complexity of
real projective space $\tc(\RP^n)$. 
This problem is
much harder than finding the topological complexity of the complex projective space (Corollary \ref{complex1}).
We will see below that the problem of finding the number $\tc(\RP^n)$ 
is equivalent to the problem of finding the smallest 
$k$ such that $\RP^n$ can be immersed into the Euclidean space $\R^{k}$.

We will begin this section by proving a general result relating the topological complexity of a space with the
Schwarz genus of a covering. 

\begin{theorem}\label{thm11} Let $X$ be a finite connected polyhedron 
and let $p:
\tilde X\to X$ be a regular covering map with the group of
covering transformations $G$. Let $\tilde X\times_G \tilde X$
be obtained from the product $\tilde X\times \tilde X$ by factorizing
with respect to the diagonal action of $G$. Then the topological
complexity $\tc(X)$ of space $X$ is greater than or equal to the
Schwarz genus of the covering
\begin{eqnarray} q: \tilde
X\times_G \tilde X\to X\times X.\label{covering}
\end{eqnarray}
\end{theorem}

\begin{proof} We will use the notation from Section \ref{section1}. In particular,
$\pi: PX \to X\times X$ will denote the canonical fibration of the space of paths
$\pi(\gamma)=(\gamma(0), \gamma(1))$, where $\gamma\in PX$, $\gamma: [0,1]\to X$. 

Consider the following commutative diagram
\begin{eqnarray*}
\begin{array}{rcl}
PX & \stackrel f\longrightarrow & \tilde X\times_G\tilde X \\ \\
& \pi\, \, \, \searrow\quad   \quad \swarrow \, \, \, q &\\ \\
& X\times X&
\end{array}
\end{eqnarray*}
where the map $f: PX \to \tilde X\times_G \tilde X$ is defined  as follows:
given a continuous path $\gamma: [0,1]\to X$, let $\tilde \gamma:[0,1]\to
\tilde X$ be any lift of $\gamma$, and we set
$$f(\gamma) =(\tilde \gamma(0), \tilde\gamma(1))\in \tilde X\times_G \tilde X.$$ 
The lift $\tilde \gamma$ of $\gamma$ depends on the choice of the initial
point $\tilde \gamma(0)\in \tilde X$ but nevertheless the map $f$ is well defined and continuous.
If $U$ is an open subset of $X\times X$ and 
$s: U\to PX$ is a continuous section of the fibration $\pi$ over $U$
then $f\circ s$ is a continuous section of $q$
over $U$. If there exits an open covering $U_1\cup U_2\cup \dots \cup U_k$
of $X\times X$ with a continuous section $s_i$ of $\pi$ over each open 
set $U_i$, then $f\circ s_i$ is a continuous section of $q$ over $U_i$
and we see that the Schwarz genus of the covering $q$ is at most $k$. 
\proofend
\end{proof}

\begin{remark}
{\rm In general (\ref{covering}) is not a regular covering;
it is regular if and only if the group $G$ is Abelian.}
\end{remark}

Next we are going to apply Theorem \ref{thm11} to the case $X=\RP^n$:

\begin{corollary}\label{new1} 
The number $\tc(\RP^n)$ is greater than or equal to the Schwarz genus of the two-fold 
covering $S^n\times_{\Z_2}S^n\to \RP^n\times \RP^n$.
\end{corollary}

We will present this result in a different form. 

If $n$ is fixed
we will always denote by 
$\xi$ the canonical real line bundle over $\RP^n$. The exterior tensor product
$\xi\otimes\xi$ is a real
line bundle over $\RP^n\times\RP^n$. Its first Stiefel - Whitney class is 
$$w_1(\xi\otimes \xi) =\alpha\times 1+ 1\times \alpha\in H^1(\RP^n\times \RP^n;\Z_2),$$ 
where $\alpha\in H^1(\RP^n;\Z_2)$ is the generator. This last condition determines uniquely the bundle $\xi\otimes\xi$.

\begin{corollary}\label{thm2} 
The topological complexity
$\tc(\RP^n)$ is not less than the minimal $k$ such that
the Whitney sum $k(\xi\otimes \xi)$ of $k$ copies of $\xi\otimes
\xi$ admits a nowhere vanishing section.
\end{corollary}

\begin{proof} By Theorem \ref{thm11} $\tc(\RP^n)$ is not less than 
the Schwarz genus of the unit sphere bundle $q$ of
$\xi\otimes \xi$. By Theorem of Schwarz \cite{Sch}, the latter
coincides with the smallest $k$ such that the $k$-fold fiberwise join $q\ast
q \ast \dots \ast q$ admits a section. But, clearly, the $k$-fold
fiberwise join $q\ast q \ast \dots \ast q$ coincides with the unit sphere
bundle of the Whitney sum $k(\xi\otimes \xi)$. This implies our statement. \proofend
\end{proof}

We will see below that the lower bound given by Corollaries \ref{new1} and \ref{thm2} are in fact precise.

The left inequality in (\ref{tab}) gives 
\begin{eqnarray}
\tc(\RP^n)\geq n+1,\label{plus}
\end{eqnarray}
since the Lusternik - Schnirelman category of $\RP^n$ is $n+1$. 
We will show later in this paper that an equality holds in (\ref{plus}) if and only if $n=1, 3, 7$.

Applying Theorem 7 from \cite{F1} we obtain:

\begin{theorem}\label{thm1}
If $n\geq 2^{r-1}$ then $\tc(\RP^n)\geq 2^r$.
\end{theorem}

\begin{proof} Let $\alpha \in H^1(\RP^n;\Z_2)$ be the generator.
Then $1\otimes \alpha +\alpha \otimes 1$ is a zero-divisor (see \S \ref{section1}), and we consider its power
\begin{eqnarray}
(1\otimes \alpha +\alpha \otimes 1)^{2^r -1}. \label{bin1}
\end{eqnarray}
The binomial expansion of this class contains the term
\begin{eqnarray}
\binom {2^r-1} n \alpha^{k} \otimes \alpha^n. \label{bin2}
\end{eqnarray}
where $k=2^r -1-n <n$.
It is well known that 
the binomial coefficients $\binom {2^r-1} i$ are odd
for all $i$. Hence (\ref{bin2}) is a non-zero term,
and so (\ref{bin1}) does not vanish either. Applying Theorem 7 from \cite{F1}, 
one finds that the topological complexity of
$\RP^n$ is not less than $2^r$. \proofend
\end{proof}

\begin{remark}{\rm 
A different proof of Theorem \ref{thm1} can be based on Corollary \ref{thm2}.
Let $k_n$ denote the least $k$ such that the top
Stiefel - Whitney class $w_k$ of $k(\xi\otimes \xi)$ vanishes.
Then $\tc(\RP^n) \geq k_n$ (since the top Stiefel - Whitney
class of a vector bundle having a section, vanishes). We want to
show that $k_n \geq 2^r$ for $n \geq 2^{r-1}$. By
Cartan's formula 
$w_k(k(\xi\otimes\xi))= (\alpha_1+\alpha_2)^k,$
where $\alpha_1, \alpha_2\in H^1(\RP^n\times \RP^n;\Z_2)$ denote
the standard 1-dimensional generators. If $k=2^r -1$, then all
binomial coefficients $\binom k  i$ are odd and since $2^r-1\leq
2n$, the Stiefel - Whitney class $w_k(k(\xi\otimes\xi))$ contains
a nontrivial term ${\binom k  i }\alpha_1^i\alpha_2^{k-i}$. Hence $k_n\geq 2^r$ and the result follows.}
\end{remark}

\section{Nonsingular maps and axial maps}

In this section we recall some notions and 
basic results concerning non-singular maps $\R^{n+1}\times \R^{n+1}\to \R^{k+1}$ and 
axial maps $\RP^n\times \RP^n\to \RP^k$. 
These maps appear in the mathematical literature in relation with the immersion problem
$\RP^n\to \R^k$, see below. The material of this section (except maybe Lemma \ref{additional})
should be considered as well-known, although we were unable to find a proper reference.

\begin{definition}\label{def1} A continuous map 
\begin{eqnarray}\label{nonsingular}
f:\R^{n}\times\R^{n}\to\R^{k}
\end{eqnarray}
is called nonsingular if it has the following two 
properties:
\begin{itemize}
\item[{\rm (a)}] $f(\lambda u, \mu v)= \lambda \mu f(u,v)$ for all $u,v\in \R^n$, 
$\lambda, \mu\in \R$, and 
\item[{\rm (b)}] $f(u,v)=0$ implies that either $u=0$, or $v=0$.
\end{itemize}
\end{definition}

Our definition of a nonsingular map is not quite standard: we do not require $f$ to be bilinear. 
Bilinear nonsingular maps give immersions of projective spaces into the Euclidean space,
see \cite{Gi}. Constructions of bilinear nonsingular maps are given in \cite{La1} and \cite{Mi}. 
In \cite{La2}, K.Y. Lam considers maps $f(u,v)$ with property (b), 
which are linear only with respect to $v$ and satisfy a
weaker property $f(-u,v)=-f(u,v)$ with respect to $u$; he calls such maps skew-linear.

We will show below (see proof of Proposition \ref{main2}) 
that the nonsingular maps in the sense of Definition \ref{def1} provide a convenient 
tool for constructing explicit motion planning algorithms in
projective spaces.

\begin{lemma}\label{borsukulam}
There are no nonsingular maps $f: \R^n\times \R^n \to \R^k$ with $k<n$.
\end{lemma}
\begin{proof} We may apply the Borsuk - Ulam Theorem to the map $u\mapsto f(u, v)$, where $v\not=0$ is fixed and
where $u$ varies on the unit sphere $S^{n-1}\subset \R^n$. By the Borsuk - Ulam Theorem, 
$f(u,v)=f(-u,v)$ for some $u\in S^{n-1}$, but the latter also is $-f(u,v)$, and thus $f(u,v)=0$.
This gives a contradiction with the nonsingularity property.
\proofend
\end{proof}

\begin{lemma}\label{dickson}
For $n=1, 2, 4$ or $8$, there exists a nonsingular map $f: \R^n \times \R^n\to \R^n$ with the property that
for any $u\in \R^n$, $u\not=0$, the first coordinate of $f(u,u)$ is positive. 
\end{lemma}
\begin{proof} For $n=1$ we take $f(u,v)=uv$, the usual product of real numbers. 

For $n=2$ we may take $f(u,v)=u\bar v$, the product of $u$ and the conjugate of $v$, viewed as complex numbers. 

We identify $\R^4$ with the set of quaternions $v=x_1+x_2i+x_3j+x_4k$ and, for $n=4$, define the nonsingular map $f$
by $f(u,v)=u\bar v\in \R^4$, where 
the bar denotes the quaternionic conjugation: $\bar v = x_1-x_2i-x_3j-x_4k$.

To construct $f$ for $n=8$, we identify $\R^8$ with the ring of Cayley numbers, which can be defined as follows, see \cite{Di}. 
A Cayley number can be uniquely written in the form $q+Qe$, where $q, Q$ are quaternions and $e$ is a formal symbol. 
The multiplication is defined by the formula
$$(q+Qe)\cdot (r+Re) = (qr-\bar RQ) +(Rq+Q\bar r)e.$$
Here the bar denotes the conjugation of quaternions, as above. We define $f:\R^8\times \R^8\to \R^8$ by
$$f(q+Qe, r+Re)=(q+Qe)\cdot (\bar r- Re).$$
Then $f(q+Qe, q+Qe)= q\bar q +Q\bar Q$ is real and positive, assuming that the Cayley number 
$q+Qe\not=0$ is nonzero. \proofend
\end{proof}

\begin{lemma}\label{adams}
For $n$ distinct from $1, 2, 4, 8$ there are no 
nonsingular maps $f: \R^n\times \R^n \to \R^n$. 
\end{lemma}
\begin{proof} The statement can be deduced from the Adams'
solution of the Hopf invariant one problem \cite{A}. 
Namely, suppose that $f$ as above exists, where $n>2$. Then we obtain a continuous map
$g: S^{n-1}\times S^{n-1}\to S^{n-1}$ given by 
$$g(x,y) = \frac{f(x,y)}{|f(x,y)|},\quad x, y\in S^{n-1}.$$
$g$ satisfies $g(-x,y) =-g(x,y)=g(x,-y)$. Restricting $g$ onto one factor $S^{n-1}\times \ast$ (where $\ast$ is a base point) 
gives a self map of $S^{n-1}$ which commutes with the antipodal involution and hence has an odd degree 
(this is a theorem of 
Borsuk, see \cite{AH}, pages 483 - 485).
Similarly, the degree of $g|_{\ast\times S^{n-1}}$ is odd. Hence, the bidegree of $g$ is $(k,\ell)$, where both integers 
$k$ and $\ell$ are odd. The existence of $g$ implies the vanishing of the Whitehead product 
\begin{eqnarray}\label{iota}
[k\iota, \ell\iota]= k\ell [\iota, \iota]\, \in \pi_{2n-3}S^{n-1},
\end{eqnarray}
where $\iota\in \pi_{n-1}S^{n-1}$ is the generator. If $n$ is odd, then $[\iota, \iota]\in \pi_{2n-3}S^{n-1}$
is of infinite order, hence the Whitehead product (\ref{iota}) cannot vanish. If $n$ is even and distinct from
$1, 2, 4, 8$, then the Whitehead product $[\iota, \iota]\in \pi_{2n-3}S^{n-1}$ is nonzero (as proven by Adams)
and has order two, which again implies that
(\ref{iota}) is nonzero as $k\ell$ is odd. \proofend
\end{proof}

\begin{definition}  Let $n$ and $k$ be two positive integers with $n<k$.
A continuous map 
\begin{eqnarray}\label{axial}
g: \RP^n\times \RP^n\to\RP^k
\end{eqnarray}
is called {\it axial of type
$(n,k)$}
if its restrictions to $*\times\RP^n$ and $\RP^n\times *$ are homotopic to the
inclusions maps $\RP^n\to \RP^k$. 
\end{definition}

Here $*$ denotes a base point of $\RP^n$.
Note that for $n<k$ any continuous map $h: \RP^n \to \RP^k$ is either homotopically trivial or it is homotopic to the 
inclusion map $\RP^n \to \RP^k$. If $\alpha\in H^1(\RP^k;\Z_2)$ denotes the generator, then 
$h^\ast\alpha\in H^1(\RP^n;\Z_2)$ is either zero or equals $\alpha$. The map $h$ is homotopically trivial
if and only if $h^\ast\alpha=0$.
This shows that the property of the axial map $g$, see (\ref{axial}), can be equivalently stated by the formula:
$g^\ast\alpha = \alpha\times 1+1\times \alpha.$ This last condition fixes the homotopy type of a map $\RP^n\times \RP^n \to \RP^\infty$ and we are interested in finding the smallest $k$ such that this map can be factorized through the inclusion
$\RP^k\to \RP^\infty$. This shows that there always exists an axial map $\RP^n\times \RP^n\to \RP^{2n}$. 
In fact, we some extra effort, one shows that there always exists an axial map $\RP^n\times \RP^n\to \RP^{2n-1}$.

\begin{lemma}\label{tab2} Assume that $1<n<k$. There exists a bijection between nonsingular maps $\R^{n+1}\times\R^{n+1}\to\R^{k+1}$ (viewed up to multiplication by a nonzero scalar) and axial maps $\RP^n\times \RP^n\to\RP^k$. 
\end{lemma}

\begin{proof}
Given a nonsingular map $f:\R^{n+1}\times\R^{n+1}\to\R^{k+1}$
one defines $g: \RP^n\times \RP^n\to\RP^k$, where for $u, v\in S^n\subset \R^{n+1}$, 
the value $g(u, v)$ is the line through the origin containing the point $f(u,v)\in \R^{k+1}$. 
To show that $g$ is indeed axial, we fix $v\in S^n$ and vary only $u\in S^n$. We see that the obtained map
$\RP^n \to \RP^k$ lifts to a map $S^n\to S^k$ given by $u\mapsto f(u,v)$, and, the relation $f(-u, v)=-f(u,v)$ implies
that $\RP^n \to \RP^k$ is not null-homotopic. Similarly, using $f(u, -v)=-f(u,v)$ we find that the restriction of $g$ onto
$\ast \times \RP^n$ is not null-homotopic. 

Suppose now that we are given an axial map (\ref{axial}). Passing to the universal covers we obtain a continuous map
$\bar g: S^n \times S^n \to S^k$ (defined up to a sign). As explained above, 
the axial property translates into $\bar g (-u,v) =-\bar g(u,v) = \bar g(u, -v)$ for all
$u, v\in S^n$. Now we may define a nonsingular map $f:\R^{n+1}\times\R^{n+1}\to\R^{k+1}$ by 
$$f(u,v) = |u|\cdot |v|\cdot \bar g\left(\frac{u}{|u|}, \frac{v}{|v|}\right), \quad u, v\in \R^{k+1}-\{0\}.$$
\proofend
\end{proof}

\begin{lemma}\label{additional} 
Suppose that for a pair of integers $1<n<k$ 
there exists a non-singular map $\R^{n+1}\times \R^{n+1}\to \R^{k+1}$. Then there exists a nonsingular map 
$f: \R^{n+1}\times \R^{n+1}\to \R^{k+1}$ having the following additional property: for any $u\in \R^{n+1}$, $u\not=0$,
the first coordinate of $f(u,u)\in \R^{k+1}$ is positive. 
\end{lemma}

\begin{proof} Given a nonsingular map $\R^{n+1}\times \R^{n+1}\to \R^{k+1}$, consider the corresponding axial map
$g: \RP^n\times \RP^n\to \RP^k$, see Lemma \ref{tab2}. 
The axial property implies that the restriction of $g$ onto the diagonal $\RP^n \subset 
\RP^n\times \RP^n$ is null-homotopic. Hence we may find $g'\simeq g$ such that 
$g':\RP^n\times \RP^n\to \RP^k$ is constant along the diagonal. Now, consider the nonsingular map 
$f: \R^{n+1}\times \R^{n+1}\to \R^{k+1}$ corresponding to $g'$ by Lemma \ref{tab2}. 
We see that for all $u\in \R^{k+1}$ 
the values $f(u,u)\in \R^{k+1}$ lie on a ray emanating from the origin. By performing an orthogonal rotation
we may assume that all nonzero vectors of this ray have positive first coordinates. This proves our claim. \proofend
\end{proof}

\section{The main results}\label{mres}

\begin{theorem}\label{main1} For any $n$, the number $\tc(\RP^n)$ coincides with the smallest 
integer $k$ such that there exists a nonsingular map
$\R^{n+1}\times \R^{n+1}\to \R^{k}$. Moreover, if $n$ is distinct from $1, 3, 7$, then $\tc(\RP^n)$ coincides with the 
smallest $k$ such that $\RP^n$ admits an immersion into $\R^{k-1}$. 
\end{theorem}

We derive from Theorem \ref{main1}, Lemma \ref{adams} and (\ref{plus}) the following Corollary:

\begin{corollary}\label{cor1} For all $n$ distinct from $1, 3,7$, one has $\tc(\RP^n)>n+1$.
The equality $\tc(\RP^n)=n+1$ takes place if and only if $n$ equals $1, 3$ or $7$. 
\end{corollary}

The second statement follows from Lemma \ref{adams} and Propositions \ref{main2}, \ref{main4} below.

Here is another Corollary:

\begin{corollary}\label{better} For any $n$, one has: $\tc({\RP^n}) \leq 2n$.
If $n$ is a power of 2 then it is an equality, i.e.
\begin{eqnarray}
\tc(\RP^{2^{r-1}})=2^r. \label{equality}
\end{eqnarray}
\end{corollary}

\begin{proof}
By the Whitney theorem the projective space $\RP^n$ (for $n>1$) admits an
immersion into  $\R^{2n-1}$, and the inequality follows (for $n\not=1, 3, 7$)
from  Theorem \ref{main1}. If $n$ is $1, 3$ or $7$ the inequality also holds since then $\tc(\RP^n)=n+1\leq 2n$.
The equality (\ref{equality}) is implied by Theorem \ref{thm1}.
 \proofend
\end{proof}

\begin{corollary}
If $n\leq n'$ then $\tc(\RP^n)\leq \tc(\RP^{n'})$.
\end{corollary}

{\it The proof} of the first statement of 
Theorem \ref{main1} for $n\not= 1, 3, 7$ follows by combining Propositions \ref{main3} and 
\ref{main2} (see below) with Lemma \ref{adams} and Corollary \ref{thm2}. 
The special cases $n=1, 3, 7$ are treated as follows:
we have the inequality (\ref{plus}) and for these values of $n$ 
we may explicitly construct motion planners with 
$n+1$ local rules (see below) which gives $\tc(\RP^n)=n+1$. 

\begin{prop}\label{main3}
For $n>1$, let $k$ be an integer such that the rank $k$ vector bundle $k(\xi\otimes \xi)$ over $\RP^n\times \RP^n$ 
(see \S \ref{initial}) admits a nowhere vanishing section. Then 
there exists a nonsingular map $\R^{n+1}\times \R^{n+1}\to \R^{k}$.
\end{prop}

\begin{prop}\label{main2} 
If there exists a nonsingular map $\R^{n+1}\times \R^{n+1}\to \R^{k},$
where $n+1<k$,
then $\RP^n$ admits 
a motion planner with $k$ local rules, i.e. 
\begin{eqnarray}
\tc(\RP^n)\leq k.\label{less}
\end{eqnarray}
\end{prop}

\begin{prop}\label{main4}
For $n=1, 3, 7$ one has $\tc(\RP^n) = n+1$.
\end{prop}

The proofs of Propositions \ref{main3}, \ref{main2} and \ref{main4} are postponed until the next section.

Propositions \ref{main3} and \ref{main2} also imply:

\begin{corollary}\label{cor8} The number 
$\tc(\RP^n)$ equals the Schwarz genus of the two-fold covering $S^n\times_{\Z_2}S^n \to \RP^n\times \RP^n$.
It also coincides with the smallest $k$ such that the vector bundle 
$k(\xi\otimes \xi)$ over $\RP^n\times \RP^n$ admits a nowhere zero section.
\end{corollary}
Two statements of Corollary \ref{cor8} are clearly equivalent, see the arguments in the proof of Corollary \ref{thm2}.

The second statement of Theorem \ref{main1} is obtained by comparing its first statement with the following
Theorem of J. Adem, S. Gitler and I.M. James and 
using the one-to-one correspondence between the axial maps and the nonsinglular maps which we have described above:

\begin{theorem} {\rm (See \cite{AGJ})} 
There exists an immersion $\RP^n\to \R^k$ (where $k>n$) if and only if there exists am axial map
$\RP^n\times \RP^n\to \RP^k$. 
\end{theorem} 

R. J. Milgram \cite{Mi} constructed, 
for any odd $n$, a nonsingular map 
$$\R^{n+1}\times \R^{n+1}\to \R^{2n+1-\alpha(n) -k(n)}.$$
Here $\alpha(n)$ denotes the number of ones in the dyadic expansion of $n$, and $k(n)$ is a non-negative function
depending only on the mod (8) residue class of $n$ with $k(1)=0$, $k(3)=k(5)=1$ and $k(7)=4$.

Applying Proposition \ref{main2} we obtain the following estimate:

\begin{corollary}\label{est} For any odd $n$ one has
$$\tc(\RP^n) \leq 2n+1-\alpha(n)-k(n).$$
\qed
\end{corollary}

Using the table from \cite{Ja1} and our Theorem \ref{main1}, one finds the topological complexity of
$\RP^n$ for all $n\leq 23$:

\begin{tabbing}
$\tc(\RP^n)$\quad \ \ \ \=2\quad \=4\quad \=4\quad \=8\quad \=8\quad
\=8\quad \=8\quad \=16\quad \=16\quad \=17\quad
\=17\quad  \=19\quad\kill
$\quad \quad n$\>1\>2\>3\>4\>5\>6\>7\>\ 8\>\ 9\>10\>11\>12\\
$\tc(\RP^n)$\>2\>4\>4\>8\>8\>8\>8\>16\>16\>17\>17\>19\\
\end{tabbing}
\begin{tabbing}
$\tc(\RP^n)$\quad \ \ \=23\quad \=23\quad \=23\quad \=32\quad
\=32\quad \=33\quad \=33\quad \=35\quad \=39\quad
\=39\quad \=39\kill
$\quad \quad n$\>13\>14\>15\>16\>17\>18\>19\>20\>21\>22\>23\\
$\tc(\RP^n)$\>23\>23\>23\>32\>32\>33\>33\>35\>39\>39\>39\\
\end{tabbing}

\section{Proofs of Propositions \ref{main3}, \ref{main2} and \ref{main4}}

{\bf Proof of Proposition \ref{main3}}.
Note that the canonical bundle $\xi$ over $\RP^n$ can be represented as the projection
$S^n\times_{\Z_2} \R\to S^n/\Z_2$.
Here $\Z_2$ acts as the antipodal involution on $S^n$ and by $x\mapsto -x$ on $\R$. 
Hence the bundle $k(\xi\otimes \xi)$ can be represented as the projection
$$(S^n\times S^n)\times_{G}\R^k\to (S^n\times S^n)/G,$$
where $G$ denotes $\Z_2\oplus \Z_2$ acting as two antipodal
involutions on $S^n\times S^n$ with each summand $\Z_2$ acting
as $x\mapsto -x$ on $\R^k$.
Consider the following commutative diagram
\begin{eqnarray*}
\begin{array}{ccc}
S^n \times S^n \times \R^k & \underset{q_1}{\to} & (S^n\times
S^n\times \R^k)/G = E(k(\xi\otimes \xi))\\ \\
\Big\downarrow \, \, p_1 &&\Big\downarrow \, \, p_2\\ \\
S^n \times S^n & \overset{q_2}{\to} & (S^n\times S^n)/G =
\RP^n\times \RP^n.
\end{array}
\end{eqnarray*}
 Then $q_1$ and $q_2$ are regular
$G$-covers, $p_1$ is the trivial $k$-plane bundle and $p_2$ can be
identified with the bundle $k(\xi\otimes \xi)$.

Any continuous section $s: \RP^n\times \RP^n \to E(k(\xi\otimes
\xi))$ of $p_2$ determines a continuous map $s_1: S^n \times  S^n
\to E(k(\xi\times \xi))$ which is invariant under $G$, i.e.
$s_1(g(x,y)) = s_1(x,y)$ for any $g\in G$. Since $n>1$ the space
$S^n\times S^n$ is simply connected and so $s_1$ can be lifted
into the covering $q_1$, which gives a continuous map $s_2:
S^n\times S^n\to S^n\times S^n \times \R^k$. The lift $s_2$ is not
unique (there are exactly four lifts) but there is a unique lift
$\tilde s_2$ having an additional property $p_1\circ \tilde
s_2=\id$, i.e. $\tilde s_2$ is a section of the trivial bundle
$p_1$. Then it is easy to see that for any $g\in G$ one has $g\circ
\tilde s_2=\tilde s_2\circ g$, i.e. $\tilde s_2$ is equivariant.

Let $f:S^n \times S^n\to \R^k -\{0\}$ be the projection of 
$s_2$ onto $\R^k$. In other words, $s_2(x,y)
=(x, y, f(x,y))$ for $x, y\in S^n$. The result is never zero
since the section $s$ assumes no zero values. We have $f(-x,y) = -f(x,y)=f(x, -y)$ for all $x,y \in S^n$.
Hence $f$
determines a nonsingular map $g:\R^{n+1}\times \R^{n+1} \to \R^k$
given by
$$g(x, y) = f\left(\frac{x}{|x|}, \frac{y}{|y|}\right), \quad x, y\in \R^{n+1}.$$
\proofend

{\bf Proof of Proposition \ref{main2}}.
We start with the following observation. Let $\phi: \R^{n+1}\times
\R^{n+1}\to \R$ be a scalar continuous map such that $\phi(\lambda u, \mu
v) =\lambda\mu \phi(u, v)$ for all $u, v\in V$ and $\lambda,
\mu\in \R$. 
Let $U_\phi\subset \RP^n\times \RP^n$
denote the set of all pairs $(A, B)$ of lines in $\R^{n+1}$ such that $A\not=B$ and
$\phi(u, v)\not=0$ for some points $u\in A$ and $v\in B$. It is
clear that $U_\phi$ is open. We claim that there exists a
continuous motion planning strategy over $U_\phi$, i.e. there is a
continuous map $s$ defined on $U_\phi$ with values in the space of continuous paths in the projective space 
$\RP^n$ such that for any pair $(A,B)\in U_\phi$ the path $s(A,B)(t)$, $t\in [0,1]$, starts at $A$ and ends at $B$. 
We may find unit
vectors $u\in A$ and $v\in B$ such that $\phi(u,v)>0$. Such pair
$u,v$ is not unique: instead of $u,v$ we may take $-u, -v$. Note
that both pairs $u, v$ and $-u, -v$ determine the same orientation
of the plane spanned by $A,B$. The desired motion planning map $s$ consists in
rotating $A$ toward $B$ in this plane, in the positive
direction determined by the orientation.

Assume now additionally that $\phi: \R^{n+1}\times \R^{n+1}\to \R$ is {\it positive} in the following sense:
for any $u\in \R^{n+1}$, $u\not=0$, one has $\phi(u,u) > 0$. Then instead of
$U_\phi$ we may take a slightly larger set $U'_\phi\subset
\RP^n\times \RP^n$, which is defined as the
 set of all pairs of lines $(A, B)$ in $\R^{n+1}$ such that
$\phi(u, v)\not=0$ for some $u\in A$ and $v\in B$. Now all pairs of lines of the form
$(A, A)$ belong to $U'_\phi$. Then for $A\not=B$ the path
from $A$ to $B$ is defined as above (rotating $A$ toward $B$ in the plane, spanned by $A$ and $B$, in the positive
direction determined by the orientation),
and for $A=B$ we choose the constant path at
$A$. Continuity is not  violated.

A vector-valued nonsingular map $f: \R^{n+1}\times \R^{n+1}\to \R^{k}$ determines $k$ scalar maps $\phi_1, \dots, \phi_k: 
\R^{n+1}\times \R^{n+1}\to \R$ (the coordinates) 
and the described above neighborhoods $U_{\phi_i}$ cover the product
$\RP^n\times \RP^n$ minus the diagonal. Since $n+1<k$ we may use Lemma \ref{additional}.
Hence we may replace the initial nonsingular map by such an $f$
that for any $u\in \R^{n+1}$, $u\not=0$, the first coordinate $\phi_1(u,u)$ of $f(u,u)$ is positive.
The open sets $U'_{\phi_1}, U_{\phi_2}, \dots, U_{\phi_k}$ cover $\RP^n\times \RP^n$. We have 
described explicit motion planning strategies over each of these sets. Therefore $\tc(\RP^n)\leq k$. 
 \proofend

{\bf Proof of Proposition \ref{main4}.} The upper bound $\tc(\RP^n)\leq n+1$ follows similarly to the proof of Proposition 
\ref{main2} and using the nonsingular maps $\R^{n+1}\times \R^{n+1}\to \R^{n+1}$ of 
Lemma \ref{dickson}. The lower bound $\tc(\RP^n)\geq n+1$ is completely general, 
see (\ref{plus}). \proofend

\section{Motion planners and immersions}
The following Theorem \ref{thm22} is a part of Theorem \ref{main1}. However in its proof we give a 
direct construction, starting from an immersion $\RP^n\to \R^k$, 
of an open covering $\{U_i\}$ of $\RP^n \times \RP^n$ with continuous motion planning strategies $s_i$ 
over each 
open set $U_i$. 

\begin{theorem}\label{thm22}
Suppose that the projective space $\RP^n$ can be immersed into 
$\R^{k}$. Then 
$\tc(\RP^n)\leq k+1.$
\end{theorem}

\begin{proof}
Imagine $\RP^n$ being immersed into $\R^k$. Fix a frame in $\R^k$ and extend it, by parallel translation, 
to a continuous field of frames. Projecting orthogonally onto $\RP^n$, 
we find $k$ continuous 
tangent vector fields $v_1, v_2, \dots, v_{k}$ on $\RP^n$ such that the vectors $v_i(p)$
(where $i=1, 2, \dots, k$) span the tangent space $T_p(\RP^{n})$ for any $p\in \RP^n$. 

Let $U_0\subset \RP^n\times \RP^n$ denote the set of pairs of lines $(A,B)$ in $\R^{n+1}$ making an acute angle. 

A nonzero tangent vector $v$ 
to the projective space $\RP^n$ at a point $A$ (which we understand as a line in $\R^{n+1}$)
determines a line $\widehat v$ in $\R^{n+1}$, which is orthogonal to $A$, i.e. $\widehat v\perp A$. The vector
$v$ also determines an orientation of the 
two-dimensional plane spanned by the lines $A$ and $\widehat v$. 

For $i=1, 2, \dots, k$
let $U_i\subset \RP^n\times \RP^n$ denote the open set of all pairs of lines $(A,B)$ in $\R^{n+1}$ 
such that the vector $v_i(A)$ is nonzero and the line
$B$ makes an acute angle with the line $\widehat{v_i(A)}$.

The sets $U_0, U_1, \dots, U_{k}$ cover $\RP^n \times \RP^n$. Indeed, given a pair $(A,B)$, there exist
indices $1\leq i_1< \dots < i_{n}\leq k$ such that the vectors $v_{i_r}(A)$, where $r=1, \dots, n$, 
span the tangent space $T_A(\RP^n)$. Then the lines 
$$A, \, \widehat {v_{i_1}(A)}, \dots, \widehat {v_{i_n}(A)}$$ 
span
the Euclidean space $\R^{n+1}$ and therefore the line $B$ makes an acute angle with one of these lines. Hence, $(A,B)$
belongs to one of the sets $U_0, U_{i_1}, \dots, U_{i_k}$. 

Now we want to describe a continuous motion planning strategy over each set $U_i$, where $i=0,1, \dots, k$.
First we define it over $U_0$. Given a pair $(A,B)\in U_0$, we rotate $A$ towards $B$ with constant velocity 
in the two-dimensional plane
spanned by $A$ and $B$ so that $A$ sweeps the acute angle. 
This clearly defines a continuous motion planning section $s_0: U_0\to P(\RP^n)$.
Our continuous motion planning strategy $s_i: U_i \to P(\RP^n)$, where $i=1, 2, \dots, k$, 
is a composition of two motions: first we rotate line $A$ toward the line $\widehat {v_i(A)}$ in the 
in the 2-dimensional 
plane spanned by $A$ and $\widehat{v_i(A)}$ in the direction determined by the orientation of this plane (see above).
On the second step we
rotate the line $\widehat {v_i(A)}$ towards $B$ along the acute angle similarly to the action of $s_0$. \proofend
\end{proof}

The inverse statement of Theorem \ref{main1}, 
allowing to construct an immersion $\RP^n \to \R^{k-1}$ starting from a motion planner for $\RP^n$,
is not very explicit; it is 
based on a long chain of constructions: Theorem \ref{thm11}, Corollary \ref{thm2}, Proposition \ref{main3}, 
and then Theorem of J. Adem, S. Gitler, I.M. James \cite{AGJ}.

\section{Motion planners in $\RP^n$ with $n\leq 7$}

Let us describe explicitly a motion planner in $\RP^2$ which may be used to solve the task of moving a line through the origin in $\R^3$.

We may view $\R^3$ as embedded into $\R^4$ (the set of quaternions) via the map 
$(x_1, x_2, x_3) \mapsto x_1+x_2i+x_3j$, where $i,j,k\in \R^4$ are the imaginary units.
Restricting the nonsingular map $\R^4\times \R^4\to \R^4$ onto $\R^3\subset \R^4$ we obtain a nonsingular map
$f:\R^3\times \R^3 \to \R^4$. Explicitly it is given by the formula
$$f(x,y) =\langle x, y\rangle -
\left|\begin{matrix}
x_1 & x_2\\
y_1 & y_2
\end{matrix}
\right|i 
-
\left|\begin{matrix}
x_1 & x_3\\
y_1 & y_3
\end{matrix}
\right|j
-
\left|\begin{matrix}
x_2 & x_3\\
y_2 & y_3
\end{matrix}
\right|k,
$$
where $x=(x_1, x_2, x_3)$ and $y=(y_1, y_2, y_3)$ and $\langle x, y\rangle$ denotes the scalar product of $x$ and $y$.
Repeating the construction given in the proof of Proposition \ref{main2} we obtain
4 open subsets $U_1, U_2, U_3, U_4$ of $\RP^2\times \RP^2$, covering 
the product $\RP^2\times \RP^2$. 
Each $U_\alpha$, where $\alpha=1, 2, 3, 4$, 
corresponds to the scalar pairing $\phi_\alpha$ obtained from $f$ by considering only the $\alpha$-th coordinate:
$$f(x,y) =\phi_1(x,y)+\phi_2(x,y)i+\phi_3(x,y)j+\phi_4(x,y)k$$
The set $U_1$ consists of the pairs of lines in $\R^3$ making an acute angle.
The set $U_2$ consists of pairs of lines in $\R^3$ such that their projections onto 
the $x_1x_2$-plane span this plane. $U_3$ and $U_4$ are defined similarly with the 
$x_1x_2$-plane replaced by the $x_2x_3$-plane and $x_1x_3$-plane, correspondingly. 

Each functionals $\phi_\alpha$ defines a continuous motion planning strategy over the 
set $U_\alpha$, see proof of Proposition \ref{main2}. For example, the motion planning strategy 
over $U_1$ is obvious: if lines $A$ and $B$ make an acute angle we rotate $A$ towards $B$ 
in the 2-plane spanned by $A$ and $B$, so that $A$ sweeps the acute angle. If $A=B$ then $A$ stays fixed.

A continuous motion planning strategy over the set $U_2$ can be described as follows. 
Let $A$ and $B$ be two lines in $\R^3$ so that their projections onto the $x_1x_2$-plane
span this plane. Fix an orientation of the $x_1x_2$-plane. For any pair $(A,B)\in U_2$ we obtain an
orientation of the 2-plane spanned by $A$ and $B$. Then we rotate $A$ towards $B$ in this 2-plane 
in the direction of the orientation.

Over the sets $U_3$ and $U_4$ we act similarly.

This example illustrates how one may use the nonsingular maps of Lemma \ref{dickson}
to construct explicitly motion planners in projective spaces $\RP^n$ with $n\leq 7$.
We skip the details. 

Notice that $\RP^3$ is homeomorphic to the Lie group ${\rm {SO}}(3)$, the configuration space of the 
rigid body, fixed at a point in 3-space. The topological complexity of this problem was computed earlier
in \S 8 of \cite{F2}. 

The immersion theory of projective spaces 
provides a rich variety of quite sophisticated 
nonsingular maps (see, for example, \cite{La1}, \cite{Mi}) which may be used to 
build, similarly to the above construction, 
explicit motion planners on higher dimensional projective spaces $\RP^n$.

\bibliographystyle{amsalpha}

\end{document}